\documentclass{article}[12pt]

\usepackage{graphicx,tikz,color}
\usepackage{amsmath,amsfonts,amssymb,latexsym}

\newtheorem{theorem}{Theorem}[section]

\newcommand{\proof}{\noindent{\bf Proof.\ }}
\newcommand{\qed}{\hfill $\square$ \bigskip}

\tikzstyle{vertex}=[circle, draw, inner sep=0pt, minimum size=6pt]

\begin{document}

\title{$M$-Polynomial Revisited: Bethe Cacti and an Extension of Gutman's Approach}

\author{Emeric Deutsch$^{a}$\thanks{Email: \texttt{EmericDeutsch@msn.com}}
\and Sandi Klav\v zar$^{b,c,d}$\thanks{Email: \texttt{sandi.klavzar@fmf.uni-lj.si}} 
}

\maketitle

\begin{center}
$^a$ Polytechnic Institute of New York University, United States\\
\medskip

$^d$ Faculty of Mathematics and Physics, University of Ljubljana, Slovenia
\medskip

$^b$ Faculty of Natural Sciences and Mathematics, University of Maribor, Slovenia\\
\medskip

$^c$ Institute of Mathematics, Physics and Mechanics, Ljubljana, Slovenia
\end{center}

\begin{abstract}
The $M$-polynomial of a graph $G$ is defined as $\sum_{i\le j} m_{i,j}(G)x^iy^j$, where $m_{i,j}(G)$, $i,j\ge 1$, is the number of edges $uv$ of $G$ such that $\{d_v(G), d_u(G)\} = \{i,j\}$. Knowing the $M$-polynomial, formulas for bond incident degree indices (an important subclass of degree-based topological indices) can be obtained by means of specific operators defined on differentiable functions in two variables. This is illustrated on three infinite families of Bethe cacti. Gutman's approach for the computation of the coefficients of the $M$-polynomial is also recalled and an extension of it is given. This extension is used to determine the $M$-polynomial of a two-parameter infinite family of lattice graphs. 
\end{abstract}


\noindent
{\bf AMS Subj.\ Class.\ (2010)}: 05C07, 05C31, 92E10.

\baselineskip16pt

\section{Introduction}
\label{sec:intro}

A large part of chemical graph theory investigates topological indices (in other words, graph invariants) which are aimed to be chemically relevant. Among these topological indices, degree-based ones, such as different variants of the Randi\'c index and the Zagreb index, play a central role. For a general and uniform discussion on the degree-based topological indices see the survey~\cite{gutman-2013a}. For selected recent investigations of  (variants of) the Randi\'c index see~\cite{chen-2017, das-2017, ma-2018, milivojevic-2017} and for (variants of) the Zagreb index we refer to~\cite{an-2018, vukicevic-2017, wang-2018}. We also refer to~\cite{gutman-2013b}, where these indices are correlated with physico-chemical parameters of octane isomers. For a selection of recent papers that compute degree-based topological indices see~\cite{ali-2018, kang-2018, kwun-2017, liu-2017, munir-2016a, munir-2016b, munir-2017, rezaei-2017}; we especially emphasize the approach to degree-based topological indices of hexagonal nanotubes in~\cite{vetrik-2017+}.

In order to simplify the computation of the bond incident degree indices, which form an utmost important subclass of degree-based topological indices  (to be defined in Section~\ref{sec:prelim}), and to stop the production of papers that, for a given family of graphs computes a given topological index from skratch, the $M$-polynomial was introduced in~\cite{deutsch-2015}. (For a related approach using the degree sequence polynomial for generalized Zagreb indices, see~\cite{doslic-2017}.) In~\cite{deutsch-2015} it was proved that the computation of several degree-based topological indices becomes a routine task, provided that the corresponding $M$-polynomial is known. More precisely, the problem more or less reduces to the one of determining the number $m_{i,j}$ of the edges of a graph whose endpoints are of degrees $i$ and $j$. Hence a particular purpose of this  paper is to point out to future authors that 
\begin{itemize}
\item[(i)] the expressions for the $m_{i,j}$s should be derived (or explained) and that 
\item[(ii)] not much space should be taken up by the computation of the topological indices; they follow easily either by elementary algebra from the $m_{i,j}$s or by elementary calculus from the $M$-polynomial.  
\end{itemize}
Numerous very recent papers that compute degree-based topological indices do not satisfy these natural requirements. Therefore, in this paper we further explain the approach and demonstrate its power on three families of Bethe cacti from~\cite{bala-1990}. These families have been selected in particular because the determination of the $M$-polynomial (equivalently of the corresponding $m_{i,j}$s) is not that straightforward as it is in several earlier cases. 

The rest of the paper is organized as follows. In the next section we formally introduce the $M$-polynomial and recall how it can be applied to bond incident degree topological indices. In Section~\ref{sec:Bethe} we introduce three families of Bethe cacti, give their recursive definitions, and based on them determine the $M$-polynomial in all of the cases. In Section~\ref{sec:Bethe-indices} we combine the results from the previous two sections to give closed formulas for several degree-based topological indices of the considered Bethe cacti. In the concluding section we recall Gutman's approach for the computation of the coefficients of the $M$-polynomial.  We extend this approach by adjoining Euler's formula to the original six equalities. We use this extended approach to determine the $M$-polynomial of a two-parameter infinite family of lattice graphs, consisting of $5$-, $6$-, and $8$-gonal faces.

We do not give basic definitions of graph theory here; the reader can consult the book~\cite{west-2001}. 

\section{Preliminaries}
\label{sec:prelim}

Let $G = (V(G), E(G))$ be a graph and let $m_{i,j}(G)$, $i,j\ge 1$, be the number of edges $uv$ of $G$ such that $\{d_v(G), d_u(G)\} = \{i,j\}$, where $d_v(G)$ (or $d_v$ for short) is the degree of the vertex $v$ in $G$. (It seems that the variables $m_{i,j}$ were introduced for the first time in~\cite{gutman-2002}.) For instance, if $G$ is $k$-regular, then $m_{k,k} = |E(G)|$, while $m_{i,j}(G) = 0$ as soon as $i\ne k$ or $j\ne k$. The {\em $M$-polynomial} of $G$ is the two variable polynomial defined as 
$$\sum_{i\le j} m_{i,j}(G)x^iy^j\,.$$ 
The role of this polynomial for degree-based indices is similar to the role of the Hosoya polynomial~\cite{hosoya-1988} (see also~\cite{deutsch-2013, eliasi-2013, lin-2013, tratnik-2017}) for distance-based invariants. 

A {\em degree-based topological index} $I$ of a graph $G$ is an arbitrary graph invariant that is defined as a function of the degrees of the vertices of $G$. In many important cases, $I$ is of the form  
\begin{equation}
\label{eq:degree-index}
I(G) = \sum_{e=uv\in E}f(d_u, d_v)\,,
\end{equation}
where $f=f(x,y)$ is a function to be suitable for chemical applications~\cite{gutman-2013a, hollas-2005}. The degree-based topological indices $I$ that are of the form~\eqref{eq:degree-index} were named {\em bond incident degree indices} in~\cite{vukicevic-2011}; we follow this terminology here. We will also abbreviate bond incident degree index to {\em BID index}. For instance, the generalized Randi\'c index $R_\alpha(G)$, $\alpha\ne 0$, is a BID index because it is obtained by selecting $f(x,y) = (xy)^\alpha$~\cite{bollobas-1998}; see Table~\ref{table1} for additional important BID indices. As examples of degree-based topological indices that are not BID indices consider the higher order Randi\'c indices. In this case the summation is taken over all paths in a graph of a given length instead over all edges as it is done in~\eqref{eq:degree-index}. 

From our point of view it is utmost important to note that~\eqref{eq:degree-index} can be rewritten as 
\begin{equation}
\label{eq:degree-index-2}
I(G) = \sum_{i\le j}m_{i,j}(G) f(i, j)\,.
\end{equation}

Consider the following operators defined on differentiable functions in two variables: 
$$
\begin{array}{cccccc}
D_x(f(x,y)) & = & x \frac{\partial f(x,y)}{\partial x}, \qquad D_y(f(x,y)) & = & y \frac{\partial f(x,y)}{\partial y}, \\ \\
S_x(f(x,y)) & = & \int_{0}^x \frac{f(t,y)}{t} dt, \qquad S_y(f(x,y)) & = & \int_{0}^y \frac{f(x,t)}{t} dt, \\ \\
J(f(x,y)) & = & f(x,x), \qquad Q_\alpha(f(x,y)) & = & x^\alpha f(x,y), \alpha\ne 0. 
\end{array}
$$
Now we can recall the following key result from~\cite{deutsch-2015}.  

\begin{theorem} {\rm \cite[Theorems 2.1,2.2]{deutsch-2015}}
\label{thm:main}
Let $G$ be a graph. 
\begin{enumerate}
\item[(i)] If $\displaystyle{I(G) = \sum_{e=uv\in E}f(d_u, d_v)}$, where $f(x,y)$ is a polynomial in $x$ and $y$, then 
$$I(G) = f(D_x,D_y)(M(G;x,y))\big|_{x=y=1}\,.$$
\item[(ii)] If $\displaystyle{I(G) = \sum_{e=uv\in E}f(d_u, d_v)}$, where $f(x,y) = \sum_{i,j\in \mathbb{Z}}\alpha_{ij}x^iy^j$, then $I(G)$ can be obtained from $M(G;x,y)$ using the operators $D_x$, $D_y$, $S_x$, and $S_y$. 
\item[(iii)] If $\displaystyle{I(G) = \sum_{e=uv\in E}f(d_u, d_v)}$, where $f(x,y) = \frac{ x^r y^s}{(x+y+\alpha)^k}$, where $r,s\ge 0$, $,t\ge 1$, and $\alpha\in {\mathbb{Z}}$, then  
$$I(G) = S_x^k\,Q_\alpha\,J\,D_x^r\,D_y^s (M(G;x,y))\big|_{x=1}\,.$$ 
\end{enumerate}
\end{theorem}

Table~\ref{table1} contains applications of Theorem~\ref{thm:main} for some of the main BID indices. 

\begin{table}[ht!]
\begin{center}
\begin{tabular}{|c|c|c|}\hline
$\phantom{\Big|}$BID index$\phantom{\Big|}$ & $f(x,y)$ & derivation from $M(G;x,y)$ \\
\hline\hline
$\phantom{\Big|}$first Zagreb$\phantom{\Big|}$ & 
$x+y$ & $(D_x + D_y)(M(G;x,y))\big|_{x=y=1}$ \\
\hline
$\phantom{\Big|}$second Zagreb$\phantom{\Big|}$ & 
$xy$ & $(D_xD_y)(M(G;x,y))\big|_{x=y=1}$ \\
\hline
$\phantom{\Big|}$second modified Zagreb$\phantom{\Big|}$ & 
$\frac{1}{xy}$ & $(S_xS_y)(M(G;x,y))\big|_{x=y=1}$ \\
\hline
$\phantom{\Big|}$general Randi\'c ($\alpha\in {\mathbb{N}}$)$\phantom{\Big|}$ & 
$(xy)^\alpha$ & $(D_x^\alpha D_y^\alpha)(M(G;x,y))\big|_{x=y=1}$ \\
\hline
$\phantom{\Big|}$general Randi\'c ($\alpha\in {\mathbb{N}}$)$\phantom{\Big|}$ & 
$\frac{1}{(xy)^\alpha}$ & $(S_x^\alpha S_y^\alpha)(M(G;x,y))\big|_{x=y=1}$ \\
\hline
$\phantom{\Big|}$symmetric division index$\phantom{\Big|}$ & 
$\frac{x^2 + y^2}{xy}$ & $(D_xS_y + D_yS_x)(M(G;x,y))\big|_{x=y=1}$ \\
\hline
$\phantom{\Big|}$harmonic$\phantom{\Big|}$ & 
$\frac{2}{x+y}$ & $2\,S_x\,J\,(M(G;x,y))\big|_{x=1}$ \\
\hline
$\phantom{\Big|}$inverse sum$\phantom{\Big|}$ & 
$\frac{xy}{x+y}$ & $S_x\,J\,D_x\,D_y\,(M(G;x,y))\big|_{x=1}$ \\
\hline
$\phantom{\Big|}$augmented Zagreb$\phantom{\Big|^{X^X}}$ & $\left(\frac{xy}{x+y-2}\right)^3$ $\phantom{\Big|}$ & $S_x^3\,Q_{-2}\,J\,D_x^3\,D_y^3\,(M(G;x,y))\big|_{x=1}$ \\
\hline
\end{tabular}
\end{center}
\caption{How to compute important BID indices from the $M$-polynomial}
\label{table1}
\end{table}

\section{Families of Bethe cacti}
\label{sec:Bethe}

Balasubramanian~\cite{bala-1990} considered families $C_n$, $D_n$, and $E_n$ ($n\ge 1$) of cactus graphs. Since the recursive structure of the families $C_n$ and $E_n$ can be described using the family $D_n$, we first consider the family $D_n$. 

\subsection{Bethe cacti $D_n$}

\begin{figure}[ht!]
\begin{center}
\begin{tikzpicture}[scale=0.9,style=thick,x=1cm,y=1cm]
\def\vr{3pt}

\coordinate(x1) at (1,1.5);
\coordinate(x2) at (1.5,1);
\coordinate(x3) at (2,1.5);
\coordinate(x4) at (1.5,2);
\coordinate(x5) at (2,2.5);
\coordinate(x6) at (2.5,2);
\coordinate(x7) at (3,2.5);
\coordinate(x8) at (2.5,3);
\coordinate(x9) at (3,1.5);
\coordinate(x10) at (3.5,1);
\coordinate(x11) at (4,1.5);
\coordinate(x12) at (3.5,2);
\coordinate(x13) at (2.5,1);
\coordinate(a1) at (5,1.5);
\coordinate(a2) at (5.5,1);
\coordinate(a3) at (6,1.5);
\coordinate(a4) at (5.5,2);
\coordinate(a5) at (6.5,0.5);
\coordinate(a6) at (7,0);
\coordinate(a7) at (7.5,0.5);
\coordinate(a8) at (7,1);
\coordinate(a9) at (6.5,2.5);
\coordinate(a10) at (7,2);
\coordinate(a11) at (7.5,2.5);
\coordinate(a12) at (7,3);
\coordinate(a13) at (7.5,3.5);
\coordinate(a14) at (8,3);
\coordinate(a15) at (8.5,3.5);
\coordinate(a16) at (8,4);
\coordinate(a17) at (8.5,5);
\coordinate(a18) at (9,4.5);
\coordinate(a19) at (9.5,5);
\coordinate(a20) at (9,5.5);
\coordinate(a21) at (9.5,3.5);
\coordinate(a22) at (10,3);
\coordinate(a23) at (10.5,3.5);
\coordinate(a24) at (10,4);
\coordinate(a25) at (10.5,2.5);
\coordinate(a26) at (11,2);
\coordinate(a27) at (11.5,2.5);
\coordinate(a28) at (11,3);
\coordinate(a29) at (10.5,0.5);
\coordinate(a30) at (11,0);
\coordinate(a31) at (11.5,0.5);
\coordinate(a32) at (11,1);
\coordinate(a33) at (12,1.5);
\coordinate(a34) at (12.5,1);
\coordinate(a35) at (13,1.5);
\coordinate(a36) at (12.5,2);
\coordinate(x) at (8,1.5);
\coordinate(y) at (9,0.5);
\coordinate(z) at (10,1.5);
\coordinate(w) at (9,2.5);
\draw (x1) -- (x2) -- (x3) -- (x4) -- (x1);
\draw (x5) -- (x6) -- (x7) -- (x8) -- (x5);
\draw (x9) -- (x10) -- (x11) -- (x12) -- (x9);
\draw (x3) -- (x6) -- (x9) -- (x13) -- (x3);
\draw (a1) -- (a2) -- (a3) -- (a4) -- (a1);
\draw (a5) -- (a6) -- (a7) -- (a8) -- (a5);
\draw (a9) -- (a10) -- (a11) -- (a12) -- (a9);
\draw (a13) -- (a14) -- (a15) -- (a16) -- (a13);
\draw (a17) -- (a18) -- (a19) -- (a20) -- (a17);
\draw (a21) -- (a22) -- (a23) -- (a24) -- (a21);
\draw (a25) -- (a26) -- (a27) -- (a28) -- (a25);
\draw (a29) -- (a30) -- (a31) -- (a32) -- (a29);
\draw (a33) -- (a34) -- (a35) -- (a36) -- (a33);
\draw (x) -- (a8) -- (a3) -- (a10) -- (x);
\draw (w) -- (a15) -- (a18) -- (a21) -- (w);
\draw (z) -- (a26) -- (a33) -- (a32) -- (z);
\draw (x) -- (y) -- (z) -- (w) -- (x);
\foreach \i in {1,2,...,13}
{
\draw(x\i)[fill=white] circle(\vr);
}
\draw(x13)[fill=black] circle(\vr);
\foreach \i in {1,2,...,36}
{
\draw(a\i)[fill=white] circle(\vr);
}
\draw(x)[fill=white] circle(\vr);
\draw(y)[fill=black] circle(\vr);
\draw(z)[fill=white] circle(\vr);
\draw(w)[fill=white] circle(\vr);
\end{tikzpicture}
\end{center}
\caption{The Bethe cacti $D_2$ and $D_3$}
\label{fig:D_2-and-D_3}
\end{figure}
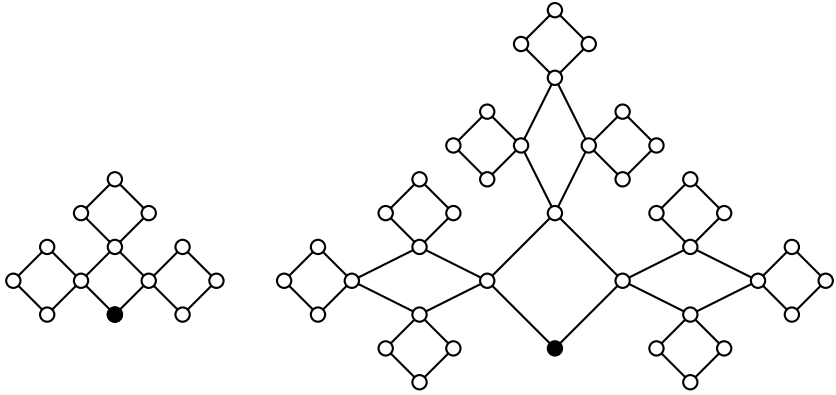

The recursive definition of the family of the Bethe cacti $D_n$, $n\ge 1$, is shown in Fig.~\ref{fig:D_n-recursive}. Here the black vertex of $D_n$ denotes the attaching vertex, where $D_n$ is attached to $D_{n+1}$ (three times). The smallest Bethe cactus $D_1$ is shown in the recursive description (Fig.~\ref{fig:D_n-recursive}), while the next two Bethe cacti $D_2$ and $D_3$ are drawn in Fig.~\ref{fig:D_2-and-D_3}. The general construction should then be clear. 

\begin{figure}[ht!]
\begin{center}
\begin{tikzpicture}[scale=0.6,style=thick,x=1cm,y=1cm]
\def\vr{3pt}

\coordinate(a1) at (0,0);
\coordinate(b1) at (1,1);
\coordinate(c1) at (0,2);
\coordinate(d1) at (-1,1);
\coordinate(a2) at (7,0);
\coordinate(b2) at (8,1);
\coordinate(c2) at (7,2);
\coordinate(d2) at (6,1);
\draw (a1) -- (b1) -- (c1) -- (d1) -- (a1);
\draw (a2) -- (b2) -- (c2) -- (d2) -- (a2);
\draw (4.6,1) ellipse (1.4cm and 1cm);
\draw (7,3.4) ellipse (1cm and 1.4cm);
\draw (9.4,1) ellipse (1.4cm and 1cm);
\foreach \i in {1,2}
{
\draw(a\i)[fill=white] circle(\vr);
\draw(b\i)[fill=white] circle(\vr);
\draw(c\i)[fill=white] circle(\vr);
\draw(d\i)[fill=white] circle(\vr);
}
\draw(a1)[fill=black] circle(\vr);
\draw(a2)[fill=black] circle(\vr);
\draw (7,-1) node {$D_n, n\ge 2$};
\draw (0,-1) node {$D_1$};
\draw (4.5,1) node {$D_{n-1}$};
\draw (9.5,1) node {$D_{n-1}$};
\draw (7,3.5) node {$D_{n-1}$};
\end{tikzpicture}
\end{center}
\caption{Recursive definition of the Bethe cacti $D_n$}
\label{fig:D_n-recursive}
\end{figure}
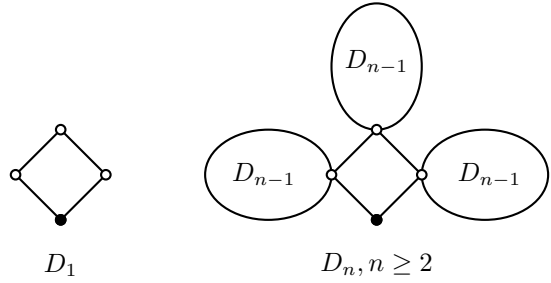

\begin{theorem}
\label{thm:M-poly-of-Dn}
$M(D_1; x,y) = 4x^2y^2$ and if $n\ge 2$, then 
$$M(D_n; x,y) = 2\cdot 3^{n-1}x^2y^2 + 2(3^{n-1}+1)x^2y^4 + 2(3^{n-1} - 2)x^4y^4\,.$$
\end{theorem}

\proof
Clearly, $M(D_1; x,y) = 4x^2y^2$. Assume in the rest that $n\ge 2$ and for the initial condition in the following three recurrences consider $D_2$ from Fig.~\ref{fig:D_2-and-D_3}. We first infer that 
$$m_{2,2}(D_{2}) = 6,\quad m_{2,2}(D_{n}) = 3 m_{2,2}(D_{n-1}), n\ge3\,,$$
which solves as $m_{2,2}(D_{n}) = 2\cdot 3^{n-1}$. 

Note further that two $24$-edges of $D_{n-1}$ become $44$-edges in $D_n$. Consequently,
$$m_{2,4}(D_{2}) = 8,\quad m_{2,4}(D_{n}) = 3 m_{2,4}(D_{n-1}) - 6 +2, n\ge3\,,$$
which solves into $m_{2,4}(D_{n}) = 2\cdot 3^{n-1} + 2$ and 
$$m_{4,4}(D_{2}) = 2,\quad m_{4,4}(D_{n}) = 3 m_{4,4}(D_{n-1}) + 6 +2, n\ge3\,,$$
which in turn solves into $m_{4,4}(D_{n}) = 2\cdot 3^{n-1} - 4$. Putting together the three solutions of the recurrences, the result follows. 
\qed

\subsection{Bethe cacti $C_n$}

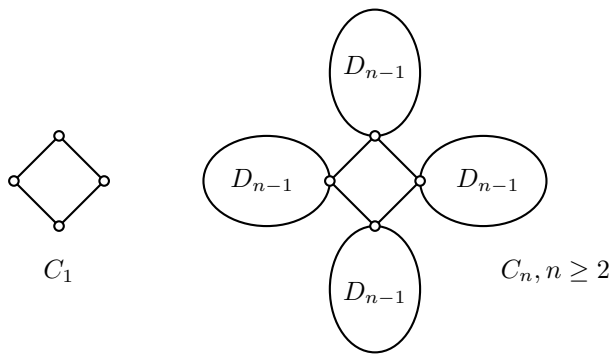
\begin{figure}[ht!]
\begin{center}
\begin{tikzpicture}[scale=0.6,style=thick,x=1cm,y=1cm]
\def\vr{3pt}

\coordinate(a1) at (0,0);
\coordinate(b1) at (1,1);
\coordinate(c1) at (0,2);
\coordinate(d1) at (-1,1);
\coordinate(a2) at (7,0);
\coordinate(b2) at (8,1);
\coordinate(c2) at (7,2);
\coordinate(d2) at (6,1);
\draw (a1) -- (b1) -- (c1) -- (d1) -- (a1);
\draw (a2) -- (b2) -- (c2) -- (d2) -- (a2);
\draw (4.6,1) ellipse (1.4cm and 1cm);
\draw (7,3.4) ellipse (1cm and 1.4cm);
\draw (9.4,1) ellipse (1.4cm and 1cm);
\draw (7,-1.4) ellipse (1cm and 1.4cm);
\foreach \i in {1,2}
{
\draw(a\i)[fill=white] circle(\vr);
\draw(b\i)[fill=white] circle(\vr);
\draw(c\i)[fill=white] circle(\vr);
\draw(d\i)[fill=white] circle(\vr);
}
\draw (11,-1) node {$C_n, n\ge 2$};
\draw (0,-1) node {$C_1$};
\draw (4.5,1) node {$D_{n-1}$};
\draw (9.5,1) node {$D_{n-1}$};
\draw (7,3.5) node {$D_{n-1}$};
\draw (7,-1.5) node {$D_{n-1}$};
\end{tikzpicture}
\end{center}
\caption{Recursive definition of the Bethe cacti $C_n$}
\label{fig:C_n-recursive}
\end{figure}

\begin{figure}[ht!]
\begin{center}
\begin{tikzpicture}[scale=0.9,style=thick,x=1cm,y=1cm]
\def\vr{3pt}

\coordinate(x1) at (1,1.5);
\coordinate(x2) at (1.5,1);
\coordinate(x3) at (2,1.5);
\coordinate(x4) at (1.5,2);
\coordinate(x5) at (2,2.5);
\coordinate(x6) at (2.5,2);
\coordinate(x7) at (3,2.5);
\coordinate(x8) at (2.5,3);
\coordinate(x9) at (3,1.5);
\coordinate(x10) at (3.5,1);
\coordinate(x11) at (4,1.5);
\coordinate(x12) at (3.5,2);
\coordinate(x13) at (2.5,1);
\coordinate(x14) at (3,0.5);
\coordinate(x15) at (2.5,0);
\coordinate(x16) at (2,0.5);

\coordinate(a1) at (5,1.5);
\coordinate(a2) at (5.5,1);
\coordinate(a3) at (6,1.5);
\coordinate(a4) at (5.5,2);
\coordinate(a5) at (6.5,0.5);
\coordinate(a6) at (7,0);
\coordinate(a7) at (7.5,0.5);
\coordinate(a8) at (7,1);
\coordinate(a9) at (6.5,2.5);
\coordinate(a10) at (7,2);
\coordinate(a11) at (7.5,2.5);
\coordinate(a12) at (7,3);
\coordinate(a13) at (7.5,3.5);
\coordinate(a14) at (8,3);
\coordinate(a15) at (8.5,3.5);
\coordinate(a16) at (8,4);
\coordinate(a17) at (8.5,5);
\coordinate(a18) at (9,4.5);
\coordinate(a19) at (9.5,5);
\coordinate(a20) at (9,5.5);
\coordinate(a21) at (9.5,3.5);
\coordinate(a22) at (10,3);
\coordinate(a23) at (10.5,3.5);
\coordinate(a24) at (10,4);
\coordinate(a25) at (10.5,2.5);
\coordinate(a26) at (11,2);
\coordinate(a27) at (11.5,2.5);
\coordinate(a28) at (11,3);
\coordinate(a29) at (10.5,0.5);
\coordinate(a30) at (11,0);
\coordinate(a31) at (11.5,0.5);
\coordinate(a32) at (11,1);
\coordinate(a33) at (12,1.5);
\coordinate(a34) at (12.5,1);
\coordinate(a35) at (13,1.5);
\coordinate(a36) at (12.5,2);
\coordinate(a37) at (7.5,-0.5);
\coordinate(a38) at (8,0);
\coordinate(a39) at (8.5,-0.5);
\coordinate(a40) at (8,-1);
\coordinate(a41) at (9.5,-0.5);
\coordinate(a42) at (10,0);
\coordinate(a43) at (10.5,-0.5);
\coordinate(a44) at (10,-1);
\coordinate(a45) at (8.5,-2);
\coordinate(a46) at (9,-1.5);
\coordinate(a47) at (9.5,-2);
\coordinate(a48) at (9,-2.5);
\coordinate(x) at (8,1.5);
\coordinate(y) at (9,0.5);
\coordinate(z) at (10,1.5);
\coordinate(w) at (9,2.5);
\draw (x1) -- (x2) -- (x3) -- (x4) -- (x1);
\draw (x5) -- (x6) -- (x7) -- (x8) -- (x5);
\draw (x9) -- (x10) -- (x11) -- (x12) -- (x9);
\draw (x3) -- (x6) -- (x9) -- (x13) -- (x3);
\draw (x13) -- (x14) -- (x15) -- (x16) -- (x13);
\draw (a1) -- (a2) -- (a3) -- (a4) -- (a1);
\draw (a5) -- (a6) -- (a7) -- (a8) -- (a5);
\draw (a9) -- (a10) -- (a11) -- (a12) -- (a9);
\draw (a13) -- (a14) -- (a15) -- (a16) -- (a13);
\draw (a17) -- (a18) -- (a19) -- (a20) -- (a17);
\draw (a21) -- (a22) -- (a23) -- (a24) -- (a21);
\draw (a25) -- (a26) -- (a27) -- (a28) -- (a25);
\draw (a29) -- (a30) -- (a31) -- (a32) -- (a29);
\draw (a33) -- (a34) -- (a35) -- (a36) -- (a33);
\draw (a37) -- (a38) -- (a39) -- (a40) -- (a37);
\draw (a41) -- (a42) -- (a43) -- (a44) -- (a41);
\draw (a45) -- (a46) -- (a47) -- (a48) -- (a45);
\draw (y) -- (a41) -- (a46) -- (a39) -- (y);
\draw (x) -- (a8) -- (a3) -- (a10) -- (x);
\draw (w) -- (a15) -- (a18) -- (a21) -- (w);
\draw (z) -- (a26) -- (a33) -- (a32) -- (z);
\draw (x) -- (y) -- (z) -- (w) -- (x);
\foreach \i in {1,2,...,16}
{
\draw(x\i)[fill=white] circle(\vr);
}
\draw(x13)[fill=white] circle(\vr);
\foreach \i in {1,2,...,48}
{
\draw(a\i)[fill=white] circle(\vr);
}
\draw(x)[fill=white] circle(\vr);
\draw(y)[fill=white] circle(\vr);
\draw(z)[fill=white] circle(\vr);
\draw(w)[fill=white] circle(\vr);
\end{tikzpicture}
\end{center}
\caption{The Bethe cacti $C_2$ and $C_3$}
\label{fig:C_2-and-C_3}
\end{figure}
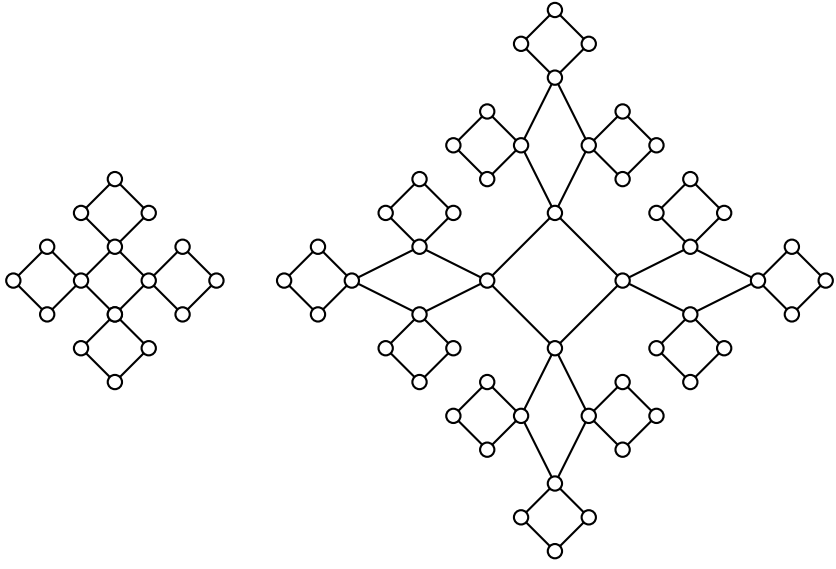

The recursive definition of the family of the Bethe cacti $C_n$, $n\ge 1$, is shown in Fig.~\ref{fig:C_n-recursive}. The vertex at which each of the four copies of $D_{n-1}$ is attached to the central $4$-cycle, respectively, is the black vertex of $D_{n-1}$ as shown in Fig.~\ref{fig:D_2-and-D_3}. The smallest Bethe cactus $C_1$ is thus the $4$-cycle graph, while the Bethe cacti $C_2$ and $C_3$ are drawn in Fig.~\ref{fig:C_2-and-C_3}. The general construction should then be clear.

\begin{theorem}
\label{thm:M-poly-of-Cn}
$M(C_1; x,y) = 4x^2y^2$ and if $n\ge 2$, then 
$$M(C_n; x,y) = 8\cdot 3^{n-2}x^2y^2 + 8\cdot 3^{n-2}x^2y^4 + 4(2\cdot 3^{n-2} - 1)x^4y^4\,.$$
\end{theorem}

\proof
Clearly, $M(C_1; x,y) = 4x^2y^2$. Assume in the rest that $n\ge 2$. Recalling from the proof of Theorem~\ref{thm:M-poly-of-Dn} that $m_{2,2}(D_{n}) = 2\cdot 3^{n-1}$, we have
$$m_{2,2}(C_{n}) = 4 m_{2,2}(D_{n-1}) = 4\cdot 2\cdot 3^{n-2} = 8\cdot 3^{n-2}\,.$$
Recalling further that $m_{2,4}(D_{n}) = 2\cdot 3^{n-1} + 2$ and $m_{4,4}(D_{n}) = 2\cdot 3^{n-1} - 4$, and observing that two $24$-edges of $D_{n-1}$ become $44$-edges in $C_n$, we get 
$$m_{2,4}(C_{n}) = 4 m_{2,4}(D_{n-1}) - 8 = 4(2\cdot 3^{n-2} + 2) - 8 = 8\cdot 3^{n-2}$$
and
$$m_{4,4}(C_{n}) = 4 m_{4,4}(D_{n-1}) + 8 + 4 = 4(2\cdot 3^{n-2} -4) + 12 = 8\cdot 3^{n-2}-4\,.$$
Hence the result. 
\qed

\subsection{Bethe cacti $E_n$}

\begin{figure}[ht!]
\begin{center}
\begin{tikzpicture}[scale=0.6,style=thick,x=1cm,y=1cm]
\def\vr{3pt}

\coordinate(a1) at (-1,1);
\coordinate(b1) at (1,1);
\coordinate(c1) at (0,2);
\coordinate(a2) at (6,1);
\coordinate(b2) at (8,1);
\coordinate(c2) at (7,2);
\draw (a1) -- (c1) -- (b1);
\draw (a2) -- (c2) -- (b2);
\draw (4.6,1) ellipse (1.4cm and 1cm);
\draw (7,3.4) ellipse (1cm and 1.4cm);
\draw (9.4,1) ellipse (1.4cm and 1cm);
\foreach \i in {1,2}
{
\draw(a\i)[fill=white] circle(\vr);
\draw(b\i)[fill=white] circle(\vr);
\draw(c\i)[fill=white] circle(\vr);
}
\draw (7,-1) node {$E_n, n\ge 2$};
\draw (0,-1) node {$E_1$};
\draw (4.5,1) node {$D_{n-1}$};
\draw (9.5,1) node {$D_{n-1}$};
\draw (7,3.5) node {$D_{n-1}$};
\end{tikzpicture}
\end{center}
\caption{Recursive definition of the Bethe cacti $E_n$}
\label{fig:E_n-recursive}
\end{figure}
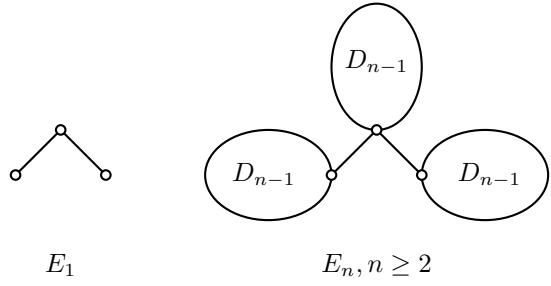

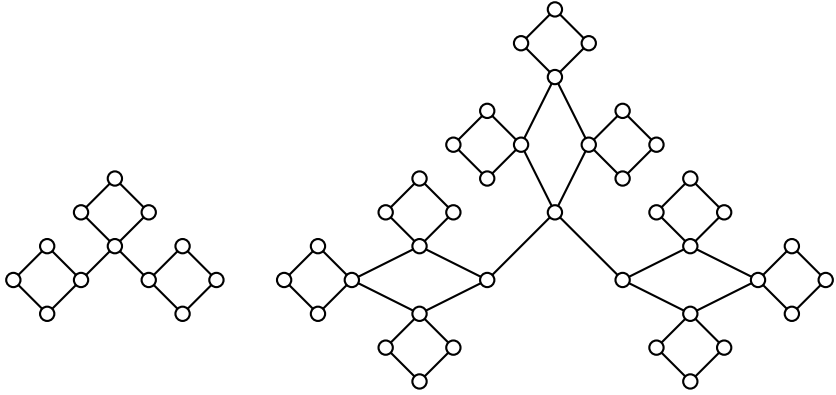
\begin{figure}[ht!]
\begin{center}
\begin{tikzpicture}[scale=0.9,style=thick,x=1cm,y=1cm]
\def\vr{3pt}

\coordinate(x1) at (1,1.5);
\coordinate(x2) at (1.5,1);
\coordinate(x3) at (2,1.5);
\coordinate(x4) at (1.5,2);
\coordinate(x5) at (2,2.5);
\coordinate(x6) at (2.5,2);
\coordinate(x7) at (3,2.5);
\coordinate(x8) at (2.5,3);
\coordinate(x9) at (3,1.5);
\coordinate(x10) at (3.5,1);
\coordinate(x11) at (4,1.5);
\coordinate(x12) at (3.5,2);
%
\coordinate(a1) at (5,1.5);
\coordinate(a2) at (5.5,1);
\coordinate(a3) at (6,1.5);
\coordinate(a4) at (5.5,2);
\coordinate(a5) at (6.5,0.5);
\coordinate(a6) at (7,0);
\coordinate(a7) at (7.5,0.5);
\coordinate(a8) at (7,1);
\coordinate(a9) at (6.5,2.5);
\coordinate(a10) at (7,2);
\coordinate(a11) at (7.5,2.5);
\coordinate(a12) at (7,3);
\coordinate(a13) at (7.5,3.5);
\coordinate(a14) at (8,3);
\coordinate(a15) at (8.5,3.5);
\coordinate(a16) at (8,4);
\coordinate(a17) at (8.5,5);
\coordinate(a18) at (9,4.5);
\coordinate(a19) at (9.5,5);
\coordinate(a20) at (9,5.5);
\coordinate(a21) at (9.5,3.5);
\coordinate(a22) at (10,3);
\coordinate(a23) at (10.5,3.5);
\coordinate(a24) at (10,4);
\coordinate(a25) at (10.5,2.5);
\coordinate(a26) at (11,2);
\coordinate(a27) at (11.5,2.5);
\coordinate(a28) at (11,3);
\coordinate(a29) at (10.5,0.5);
\coordinate(a30) at (11,0);
\coordinate(a31) at (11.5,0.5);
\coordinate(a32) at (11,1);
\coordinate(a33) at (12,1.5);
\coordinate(a34) at (12.5,1);
\coordinate(a35) at (13,1.5);
\coordinate(a36) at (12.5,2);
\coordinate(x) at (8,1.5);
\coordinate(z) at (10,1.5);
\coordinate(w) at (9,2.5);
\draw (x1) -- (x2) -- (x3) -- (x4) -- (x1);
\draw (x5) -- (x6) -- (x7) -- (x8) -- (x5);
\draw (x9) -- (x10) -- (x11) -- (x12) -- (x9);
\draw (x3) -- (x6) -- (x9);
\draw (a1) -- (a2) -- (a3) -- (a4) -- (a1);
\draw (a5) -- (a6) -- (a7) -- (a8) -- (a5);
\draw (a9) -- (a10) -- (a11) -- (a12) -- (a9);
\draw (a13) -- (a14) -- (a15) -- (a16) -- (a13);
\draw (a17) -- (a18) -- (a19) -- (a20) -- (a17);
\draw (a21) -- (a22) -- (a23) -- (a24) -- (a21);
\draw (a25) -- (a26) -- (a27) -- (a28) -- (a25);
\draw (a29) -- (a30) -- (a31) -- (a32) -- (a29);
\draw (a33) -- (a34) -- (a35) -- (a36) -- (a33);
\draw (x) -- (a8) -- (a3) -- (a10) -- (x);
\draw (w) -- (a15) -- (a18) -- (a21) -- (w);
\draw (z) -- (a26) -- (a33) -- (a32) -- (z);
\draw (z) -- (w) -- (x);
\foreach \i in {1,2,...,12}
{
\draw(x\i)[fill=white] circle(\vr);
}
\foreach \i in {1,2,...,36}
{
\draw(a\i)[fill=white] circle(\vr);
}
\draw(x)[fill=white] circle(\vr);
\draw(z)[fill=white] circle(\vr);
\draw(w)[fill=white] circle(\vr);
\end{tikzpicture}
\end{center}
\caption{The Bethe cacti $E_2$ and $E_3$}
\label{fig:E_2-and-E_3}
\end{figure}

The recursive definition of the family of the Bethe cacti $E_n$, $n\ge 1$, is shown in Fig.~\ref{fig:E_n-recursive}. Again, the vertex at which each of the three copies of $D_{n-1}$ is attached to the central path on three vertices, respectively, is the black vertex of $D_{n-1}$ as shown in Fig.~\ref{fig:D_2-and-D_3}. Thus the smallest Bethe cactus $E_1$ is the path on three vertices, the next two Bethe cacti $E_2$ and $E_3$ are drawn in Fig.~\ref{fig:E_2-and-E_3}. The general construction should then be clear.

\begin{theorem}
\label{thm:M-poly-of-En}
$M(E_1; x,y) = 2xy^2$, $M(E_2; x,y) = 6x^2y^2 + 4x^2y^3 + 2x^2y^4 + 2x^3y^4$, and if $n\ge 3$, then 
$$M(E_n; x,y) = 2\cdot 3^{n-1}x^2y^2 + 2\cdot 3^{n-1}x^2y^4 + 6x^3y^4 + (2\cdot 3^{n-1} - 10)x^4y^4\,.$$
\end{theorem}

\proof
Clearly, $M(E_1; x,y) = 2xy^2$ and $M(E_2; x,y) = 6x^2y^2 + 4x^2y^3 + 2x^2y^4 + 2x^3y^4$. Assume in the rest that $n\ge 3$. Note that two $24$-edges of the middle $D_{n-1}$ become $44$-edges in $E_n$, and that two $24$-edges of an extreme $D_{n-1}$ become $34$-edges in $E_n$. Hence, recalling again from the proof of Theorem~\ref{thm:M-poly-of-Dn} that $m_{2,2}(D_{n}) = 2\cdot 3^{n-1}$, $m_{2,4}(D_{n}) = 2\cdot 3^{n-1} + 2$, and $m_{4,4}(D_{n}) = 2\cdot 3^{n-1} - 4$, we get: 
\begin{eqnarray*}
m_{2,2}(E_{n}) & = & 3 m_{2,2}(D_{n-1}) = 3\cdot 2\cdot 3^{n-2} = 2\cdot 3^{n-1}\,,\\
m_{2,3}(E_{n}) & = & 3 m_{2,3}(D_{n-1}) = 0\,,\\
m_{2,4}(E_{n}) & = & 3 m_{2,4}(D_{n-1}) - 6 = 3\cdot 2\cdot 3^{n-2} + 6 - 6 = 2\cdot 3^{n-1}\,,\\
m_{3,3}(E_{n}) & = & 3 m_{3,3}(D_{n-1}) = 0\,,\\
m_{3,4}(E_{n}) & = & 3 m_{3,4}(D_{n-1}) + 4 + 2 = 6\,,\\
m_{4,4}(E_{n}) & = & 3 m_{4,4}(D_{n-1}) + 2 = 3\cdot 2\cdot 3^{n-2} -12 + 2 = 2\cdot 3^{n-1} - 10\,.
\end{eqnarray*}
Putting all this together, the result follows. 
\qed

\section{Topological indices of Bethe cacti}
\label{sec:Bethe-indices}

Combining Theorem~\ref{thm:M-poly-of-Dn} with the expressions from Table~\ref{table1}, routine computations yield the expressions for the selected listed topological indices of $D_n$, $C_n$, and $E_n$, $n\ge 2$. Let us demonstrate this by computing the symmetric division index of $D_n$, $n\ge 2$. From Table~\ref{table1} we know that this reduces to compute $(D_xS_y + D_yS_x)(M(G;x,y))\big|_{x=y=1}$. Now,
\begin{eqnarray*}
S_x(M(D_n;x,y)) & = & \int_0^x \frac{2\cdot 3^{n-1}t^2y^2 + 2(3^{n-1}+1)t^2y^4 + 2(3^{n-1} - 2)t^4y^4}{t} \\
& = & \frac{3^{n - 1}x^2y^2(x^2y^2 + 2(y^2 + 1))}{2} - x^4y^4 + x^2y^4\,,
\end{eqnarray*}
and hence
\begin{eqnarray}
\label{eq:first}
D_yS_x(M(D_n;x,y)) & = & y\cdot (2\cdot 3^{n - 1}x^2y(x^2y^2 + 2y^2 + 1) - 4x^2y^3(x^2 - 1))\,.
\end{eqnarray}
Similarly we compute that 
\begin{eqnarray}
\label{eq:second}
D_xS_y(M(D_n;x,y)) & = & x\cdot (3^{n - 1}xy^2(2x^2y^2 + y^2 + 2) - xy^4(4x^2 - 1))\,.
\end{eqnarray}
Summing~\eqref{eq:first} and~\eqref{eq:second} we get
\begin{eqnarray*}
(D_xS_y + D_yS_x)(M(G;x,y)) & = & 3^{n - 1}x^2y^2(4x^2y^2 + 5y^2 + 4) - x^2y^4·(8x^2 - 5)
\end{eqnarray*}
from where we conclude that 
\begin{eqnarray*}
(D_xS_y + D_yS_x)(M(G;x,y))\big|_{x=y=1} & = & 13\cdot 3^{n - 1} - 3\,.
\end{eqnarray*}
All the other entries from Table~\ref{table2} are computed along the same lines. 

\begin{table}[ht!]
\begin{center}
\begin{tabular}{|c|c|c|c|}\hline
$\phantom{\Big|}$topological index $I$ & $I(D_n)$ & $I(C_n)$ & $I(E_n)$ \\
\hline\hline
first Zagreb & 
$\phantom{\displaystyle{\frac{1}{1}}} 4\cdot3^{n + 1} - 20$ & 
$16\cdot 3^n - 32$ & 
$4\cdot 3^{n + 1} - 38$ \\
\hline
second Zagreb & 
$\phantom{\displaystyle{\frac{1}{1}}} 56\cdot3^{n - 1} - 48$ & 
$224\cdot 3^{n - 2} - 64$ & 
$56\cdot 3^{n - 1} - 88$ \\
\hline
second modified Zagreb & 
$\phantom{\displaystyle{\frac{1}{1}}} \frac{7}{8}\cdot 3^{n - 1}$ & 
$\frac{7}{2}\cdot 3^{n - 2} - \frac{1}{4}$ & 
$\frac{7}{8}\cdot 3^{n - 1} - \frac{1}{8}$ \\
\hline
symmetric division index & 
$\phantom{\displaystyle{\frac{1}{1}}} 13\cdot3^{n - 1} - 3$ & 
$52\cdot 3^{n - 2} - 8$ & 
$13\cdot 3^{n - 1} - \frac{15}{2}$ \\
\hline
harmonic & 
$\phantom{\displaystyle{\frac{1}{1}}} \frac{13}{2}\cdot 3^{n - 2} - \frac{1}{3}$ &
$26\cdot 3^{n - 3} - 1$ & 
$\frac{13}{2}\cdot 3^{n - 2} - \frac{11}{14}$ \\
\hline
inverse sum & 
$\phantom{\displaystyle{\frac{1}{1}}} 26\cdot 3^{n - 2} - \frac{16}{3}$ &
$104\cdot 3^{n - 3} - 8$ & 
$26\cdot 3^{n - 2} - \frac{68}{7}$ \\
\hline
\end{tabular}
\end{center}
\caption{Selected topological indices of Bethe cacti}
\label{table2}
\end{table}

\section{(An extension of) Gutman's approach}
\label{sec:conclude}

As already pointed out in~\cite{deutsch-2015}, an approach to determine the coefficients $m_{i,j}$ of an $M$-polynomial has been proposed by Gutman~\cite{gutman-2002} by considering corresponding linear equations. Let us briefly recall the approach here, in particular to correct a statement from~\cite[p.~99]{deutsch-2015} (see below). 

Let $G$ be a chemical graph (a graph of maximum degree at most $4$) with $n$ vertices and $m$ edges, and let $n_i$, $1\le i\le 4$, be the number of vertices of degree $i$. Clearly, $m_{1,1} = 0$ as soon as the graph has at least three vertices and is connected, while for the the other $m_{i,j}$s we have: 
\begin{eqnarray}
n_1 + n_2 + n_3 + n_4 & = & n \label{eq-3} \\
m_{1,2} + m_{1,3} + m_{1,4} & = & n_1 \label{eq-4} \\
m_{1,2} + 2m_{2,2} + m_{2,3} + m_{2,4} & = & 2n_2 \label{eq-5} \\
m_{1,3} + m_{2,3} + 2m_{3,3} + m_{3,4} & = & 3n_3 \label{eq-6} \\
m_{1,4} + m_{2,4} +m_{3,4} + 2m_{4,4} & = & 4n_4  \label{eq-7}\\
n_1 + 2n_2 + 3n_3 + 4n_4 & = & 2m\,. \label{eq-8}
\end{eqnarray}
Equations~\eqref{eq-3}-\eqref{eq-7} are linearly independent, while~\eqref{eq-8} is a consequence of~\eqref{eq-3}-\eqref{eq-7}. (In~\cite{deutsch-2015} it is said that all these equations are linearly independent.) Gutman's approach is to determine first some of the $m_{i,j}$s and then the remaining ones can be obtained from the above relations.  

We extend Gutman's approach by adjoining to Equations~\eqref{eq-3}-\eqref{eq-8} Euler's formula (cf.~\cite[p.~201]{west-2001})
\begin{equation}
\label{eq:Euler}
\sum m_{i,j} - \sum n_i = f - 2\,,
\end{equation}
usable whenever dealing with a plane graph whose number of faces $f$ can be determined.

In the rest we are going to use this extended Gutman approach to determine the $M$-polynomial of the networks $G(p,q)$, $p,q\ge 1$. In Fig.~\ref{fig:G(3,4)} the network $G(3,4)$ is drawn, from which the general definition should be clear. In particular, $G(1,1)$ consists of an $8$-gon with two $6$-gons attached at the top and two $6$-gons attached at the bottom.

\begin{figure}[ht!]
\begin{center}
\begin{tikzpicture}[scale=0.7,style=thick,x=1cm,y=1cm]
\def\vr{3pt}
\draw (1,0) -- (0,0.5) -- (0,1.5) -- (1,2) -- (0.5,3) -- (1,4) -- (0,4.5) -- (0,5.5) -- (1,6) -- (0.5,7) -- (1,8) -- (0,8.5) -- (0,9.5) -- (1,10) -- (0.5,11) -- (1,12) -- (0,12.5) -- (0,13.5) -- (1,14);
\draw (5,0) -- (4,0.5) -- (4,1.5) -- (5,2) -- (4.5,3) -- (5,4) -- (4,4.5) -- (4,5.5) -- (5,6) -- (4.5,7) -- (5,8) -- (4,8.5) -- (4,9.5) -- (5,10) -- (4.5,11) -- (5,12) -- (4,12.5) -- (4,13.5) -- (5,14);
\draw (9,0) -- (8,0.5) -- (8,1.5) -- (9,2) -- (8.5,3) -- (9,4) -- (8,4.5) -- (8,5.5) -- (9,6) -- (8.5,7) -- (9,8) -- (8,8.5) -- (8,9.5) -- (9,10) -- (8.5,11) -- (9,12) -- (8,12.5) -- (8,13.5) -- (9,14);
\draw (13,0) -- (12,0.5) -- (12,1.5) -- (13,2) -- (12.5,3) -- (13,4) -- (12,4.5) -- (12,5.5) -- (13,6) -- (12.5,7) -- (13,8) -- (12,8.5) -- (12,9.5) -- (13,10) -- (12.5,11) -- (13,12) -- (12,12.5) -- (12,13.5) -- (13,14);
\draw (1,0) -- (2,0.5) -- (2,1.5) -- (3,2) -- (3.5,3) -- (3,4) -- (2,4.5) -- (2,5.5) -- (3,6) -- (3.5,7) -- (3,8) -- (2,8.5) -- (2,9.5) -- (3,10) -- (3.5,11) -- (3,12) -- (2,12.5) -- (2,13.5) -- (1,14);
\draw (5,0) -- (6,0.5) -- (6,1.5) -- (7,2) -- (7.5,3) -- (7,4) -- (6,4.5) -- (6,5.5) -- (7,6) -- (7.5,7) -- (7,8) -- (6,8.5) -- (6,9.5) -- (7,10) -- (7.5,11) -- (7,12) -- (6,12.5) -- (6,13.5) -- (5,14);
\draw (9,0) -- (10,0.5) -- (10,1.5) -- (11,2) -- (11.5,3) -- (11,4) -- (10,4.5) -- (10,5.5) -- (11,6) -- (11.5,7) -- (11,8) -- (10,8.5) -- (10,9.5) -- (11,10) -- (11.5,11) -- (11,12) -- (10,12.5) -- (10,13.5) -- (9,14);
\draw (13,0) -- (14,0.5) -- (14,1.5) -- (15,2) -- (15.5,3) -- (15,4) -- (14,4.5) -- (14,5.5) -- (15,6) -- (15.5,7) -- (15,8) -- (14,8.5) -- (14,9.5) -- (15,10) -- (15.5,11) -- (15,12) -- (14,12.5) -- (14,13.5) -- (13,14);
\draw (0,0.5) -- (1,0) -- (2,0.5) -- (3,0) -- (4,0.5) -- (5,0) -- (6,0.5) -- (7,0) -- (8,0.5) -- (9,0) -- (10,0.5) -- (11,0) -- (12,0.5) -- (13,0) -- (14,0.5) -- (15,0) -- (16,0.5);
\draw (0,4.5) -- (1,4) -- (2,4.5) -- (3,4) -- (4,4.5) -- (5,4) -- (6,4.5) -- (7,4) -- (8,4.5) -- (9,4) -- (10,4.5) -- (11,4) -- (12,4.5) -- (13,4) -- (14,4.5) -- (15,4) -- (16,4.5);
\draw (0,8.5) -- (1,8) -- (2,8.5) -- (3,8) -- (4,8.5) -- (5,8) -- (6,8.5) -- (7,8) -- (8,8.5) -- (9,8) -- (10,8.5) -- (11,8) -- (12,8.5) -- (13,8) -- (14,8.5) -- (15,8) -- (16,8.5);
\draw (0,12.5) -- (1,12) -- (2,12.5) -- (3,12) -- (4,12.5) -- (5,12) -- (6,12.5) -- (7,12) -- (8,12.5) -- (9,12) -- (10,12.5) -- (11,12) -- (12,12.5) -- (13,12) -- (14,12.5) -- (15,12) -- (16,12.5);
\draw (0,1.5) -- (1,2) -- (2,1.5) -- (3,2) -- (4,1.5) -- (5,2) -- (6,1.5) -- (7,2) -- (8,1.5) -- (9,2) -- (10,1.5) -- (11,2) -- (12,1.5) -- (13,2) -- (14,1.5) -- (15,2) -- (16,1.5);
\draw (0,5.5) -- (1,6) -- (2,5.5) -- (3,6) -- (4,5.5) -- (5,6) -- (6,5.5) -- (7,6) -- (8,5.5) -- (9,6) -- (10,5.5) -- (11,6) -- (12,5.5) -- (13,6) -- (14,5.5) -- (15,6) -- (16,5.5);
\draw (0,9.5) -- (1,10) -- (2,9.5) -- (3,10) -- (4,9.5) -- (5,10) -- (6,9.5) -- (7,10) -- (8,9.5) -- (9,10) -- (10,9.5) -- (11,10) -- (12,9.5) -- (13,10) -- (14,9.5) -- (15,10) -- (16,9.5);
\draw (0,13.5) -- (1,14) -- (2,13.5) -- (3,14) -- (4,13.5) -- (5,14) -- (6,13.5) -- (7,14) -- (8,13.5) -- (9,14) -- (10,13.5) -- (11,14) -- (12,13.5) -- (13,14) -- (14,13.5) -- (15,14) -- (16,13.5);
\draw (3.5,3) -- (4.5,3); \draw (7.5,3) -- (8.5,3); \draw (11.5,3) -- (12.5,3); 
\draw (3.5,7) -- (4.5,7); \draw (7.5,7) -- (8.5,7); \draw (11.5,7) -- (12.5,7); 
\draw (3.5,11) -- (4.5,11); \draw (7.5,11) -- (8.5,11); \draw (11.5,11) -- (12.5,11); 
\draw (16,0.5) -- (16,1.5); 
\draw (16,4.5) -- (16,5.5); 
\draw (16,8.5) -- (16,9.5); 
\draw (16,12.5) -- (16,13.5); 
\draw (2,11) node {$1$};
\draw (2,7) node {$2$};
\draw (2,3) node {$p=3$};
\draw (6,11) node {$2$};
\draw (10,11) node {$3$};
\draw (14,11) node {$q=4$};

\end{tikzpicture}
\end{center}
\caption{The lattice $G(3,4)$}
\label{fig:G(3,4)}
\end{figure}
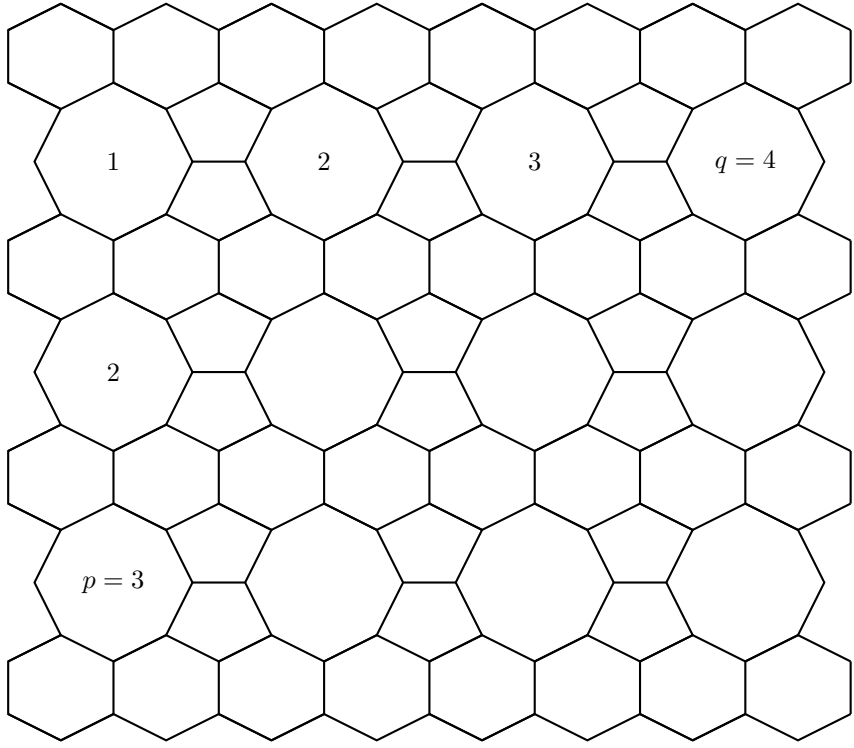

Clearly, vertices of $G(p,q)$ are of degrees $2$ and $3$, hence we need to determine $m_{2,2} = m_{2,2}(G(p,q))$, 
$m_{2,3} = m_{2,3}(G(p,q))$, and
$m_{3,3} = m_{3,3}(G(p,q))$. 

Note first that 
\begin{equation}
\label{eq:m22}
m_{2,2} = 2(p + 1) + 4 = 2p + 6\,,
\end{equation}
where $2(p + 1)$ correspond the side edges with both end-points of degree $2$, and $4$ corresponds to the corner edges (with both end-points of degree $2$). Furthermore,  $n_2 = 4q + 4(p + 1) + 2p = 6p + 4q + 4$, where $4q$ comes from the top and bottom vertices of degree $2$, the term $4(p + 1)$ comes from the sides, and the term $2p$ from the almost sides. Equation~\eqref{eq-5} in our case reduces to $2m_{2,2} + m_{2,3} = 2n_2$, from which we get
\begin{equation}
\label{eq:m23}
m_{2,3} = 8p + 8q - 4\,.
\end{equation}
Equation~\eqref{eq-6} reduces to $m_{2,3} + 2m_{3,3} = 3n_3$ and therefore, 
\begin{equation}
\label{eq:n3-m33}
3n_3 - 2m_{3,3} = 8p + 8q - 4\,.
\end{equation}
Since the number of $8$-gons of $G(p,q)$ is $pq$, the number of its $6$-gons is $2q(p+1)$, and the number of its $5$-gons is $2p(q-1)$, Equation~\eqref{eq:Euler} reduces to 
\begin{equation}
\label{eq:Euler-2}
m_{3,3} - n_3 = 5pq - 6p - 2q + 1\,.
\end{equation}
Solving~\eqref{eq:n3-m33} and~\eqref{eq:Euler-2} yields $n_3 = 10pq - 4p + 4q - 2$ and 
\begin{equation}
\label{eq:m33}
m_{3,3} = 15pq - 10p + 2q - 1\,.
\end{equation}
From Equations~\eqref{eq:m22}, \eqref{eq:m23}, and \eqref{eq:m33} we conclude that 
$$M(G(p,q); x,y) = (2p + 6)x^2y^2 + (8p + 8q - 4)x^2y^3 + (15pq - 10p + 2q - 1)x^3y^3\,.$$

\section*{Acknowledgments}

Sandi Klav\v zar acknowledges the financial support from the Slovenian Research Agency (research core funding No.\ P1-0297).

\end{document}